\documentclass[a4paper, 12pt, leqno]{amsart}

\usepackage{amsmath,amssymb,amsthm,eucal,bm}
\usepackage{graphicx}
\usepackage{cite,mathrsfs}

\DeclareFontEncoding{OT2}{}{}
\DeclareFontSubstitution{OT2}{cmr}{m}{n}
\DeclareSymbolFont{cyss}{OT2}{wncyss}{m}{n}
\DeclareMathSymbol{\sh}{\mathbin}{cyss}{`x}

\newtheorem{theorem}{Theorem}[section]

\theoremstyle{definition}

\newtheorem{assumption}[theorem]{Assumption}

\theoremstyle{remark}
\newtheorem{remark}[theorem]{Remark}

\numberwithin{equation}{section}

\begin{document}

\title[Multiple Dirichlet series with arithmetical coefficients]{A survey on the theory of multiple Dirichlet series with arithmetical coefficients
on the numerators}

\author{Kohji Matsumoto}
\address{K. Matsumoto: Graduate School of Mathematics, Nagoya University, Chikusa-\
ku, Nagoya 464-8602, Japan}
\email{kohjimat@math.nagoya-u.ac.jp}


\maketitle
\baselineskip 16pt

\begin{center}
{\it Dedicated to Professor Jonas Kubilius \\on the occasion of the 100 years 
anniversaity of his birth}
\end{center}

\begin{abstract}
We survey some recent developments in the analytic theory of multiple Dirichlet series with
arithmetical coefficients on the numerators.
\end{abstract}
\bigskip

\section{A personal recollection}

I first met Professor Jonas Kubilius at Kyoto in July 1986, when the
5th USSR-Japan Symposium on Probability was held there.     
Professor Kubilius was one of the members of
the Soviet team, and on this occasion, some Japanese number theorists organized a small
satellite meeting on probabilistic number theory with him.    At that time I was
a post-doctoral researcher, just after getting my degree of Dr.~Sci. from Rikkyo
University in March of the same
year.    On the meeting I gave a talk on the contents of my thesis, concerning the
value-distribution of the Riemann zeta-function $\zeta(s)$.    After my talk, Professor
Kubilius came close to me and said:

``Do you know the name of Antanas Laurin{\v c}ikas?''

``No.''

``He is a Lithuanian mathematician, and he wrote a lot of papers in which you are surely
interested.''

Further he mentioned that many papers of Laurin{\v c}ikas can be found in Lietuvos
Matematikos Rinkinys.   I thanked Professor Kubilius, and when I went 
back to Rikkyo University, I immediately visited the library.     Unfortunately in the
library there were only the Russian original version of the journal back numbers, which
I could not read, so at that time I gave up reading the papers of Laurin{\v c}ikas.
(Several years later I noticed the existence of the English translation of Liet. Mat.
Rink.)

In the next year, I got a job of a lecturer at Iwate University, and in 1995 I moved to
Nagoya University.   And then in 1996 I first visited Lithuania, to attend the 2nd
Palanga Conference on Analytic and Probabilistic Number Theory.    I met again 
Professor Kubilius, and found many new Lithuanian friends.
Since then, I became a regular member of the Palanga Conferences.    
In September 2011, I attended the 5th Palanga Conference; it was just one and a half month
before the death of Professor Kubilius.    When I arrived at the university villa in
Palanga, at the entrance of the villa I met Professor Kubilius.     
When he noticed me, he said

``{\c C}a va?''

I replied oui {\c ca} va, and said greetings in English.
It was my last conversation with him.

\section{Multiple Dirichlet series with or without coefficients}

One of the most favorite topics investigated by Professor Kubilius is the theory of 
arithmetical functions.
Therefore the author feels it a natural choice here to report some recent developments in
the theory of multiple Dirichlet series with arithmetical coefficients on the
numerators.

We begin, however, with the definition of multiple zeta-functions without coefficients:
\begin{align}\label{def_EZ}
&\zeta_r(s_1,\ldots,s_r)=\sum_{m_1=1}^{\infty}\cdots \sum_{m_r=1}^{\infty}
m_1^{-s_1}(m_1+m_2)^{-s_2}\\
&\qquad\qquad\times\cdots\times (m_1+\cdots+m_r)^{-s_r}, \notag
\end{align}
where $s_1,\ldots,s_r\in\mathbb{C}$.
This multiple series is sometimes called the Euler-Zagier $r$-fold zeta-function.
This series is absolutely convergent when
$$
\Re(s_{r-k+1}+\cdots+s_r)>k \qquad (1\leq k\leq r),
$$
but can be continued meromorphically to the whole space $\mathbb{C}^r$.
Now lots of analytic, algebraic, and arithmetic properties of this series \eqref{def_EZ}
and its special values are known.

It is natural to consider some generalization of \eqref{def_EZ}, which has some
coefficients on the numerators.    We introduce the following two types of
generalizations:

\noindent{\bf Type $\sh$}:
\begin{align*}
&\Phi_r^{\sh}(s_1,\ldots,s_r;a_1,\ldots,a_r)\\
&\quad=\sum_{m_1=1}^{\infty}\cdots\sum_{m_r=1}^{\infty}
\frac{a_1(m_1)a_2(m_2)\cdots a_r(m_r)}{m_1^{s_1}(m_1+m_2)^{s_2}\cdots(m_1+\cdots+m_r)^{s_r}},
\end{align*}
and

\noindent{\bf Type $*$}:
\begin{align*}
&\Phi_r^{*}(s_1,\ldots,s_r;a_1,\ldots,a_r)\\
&\quad=\sum_{m_1=1}^{\infty}\cdots\sum_{m_r=1}^{\infty}
\frac{a_1(m_1)a_2(m_1+m_2)\cdots a_r(m_1+\cdots+m_r)}{m_1^{s_1}(m_1+m_2)^{s_2}
\cdots(m_1+\cdots+m_r)^{s_r}}.
\end{align*}

Here, the notation $\sh$ and $*$ is according to the philosophy of T. Arakawa and M. Kaneko
\cite{ArKa04}.
The aim of the present article is to survey recent results related with the above two
types of multiple series.     

It is to be noted that the present survey is by no means
complete.   
In this article we mainly discuss the analytic point of view, so many important results
on special values are not mentioned.
Moreover, we do not discuss several multiple series with arithmetical numerators,
similar to our series.
For example, there are a lot of references on multiple Dirichlet series with some kind of
twisted factors on the numerators.    In the present article, however, we will only
discuss the case of Dirichlet characters later, without mentioning the details of
other articles,
e.g., P. Cassou-Nogu{\`e}s \cite{CN81} \cite{CN8287}, M. de Crisenoy \cite{dC06},
de Crisenoy and D. Essouabri \cite{dCEs08}, Essouabri and the author \cite{EsMa19}, and so on.

More generally, 
R. de la Bret{\`e}che \cite{dlBr01} studied the multiple series of the form
$$
\sum_{m_1=1}^{\infty}\cdots\sum_{m_r=1}^{\infty}
\frac{f(m_1,\ldots,m_r)}{m_1^{s_1}m_2^{s_2}\cdots m_r^{s_r}},
$$ 
while several mathematicians treated 
$$
\sum_{m_1=1}^{\infty}\cdots\sum_{m_r=1}^{\infty}
\frac{f(m_1,\ldots,m_r)}{P(m_1,\ldots,m_r)^s}
$$
where $P$ is a polynomial, or more generally
$$
\sum_{m_1=1}^{\infty}\cdots\sum_{m_r=1}^{\infty}
\frac{f(m_1,\ldots,m_r)}{P_1(m_1,\ldots,m_r)^{s_1}\cdots P_n(m_1,\ldots,m_r)^{s_n}}
$$
where $P_1,\ldots,P_n$ are polynomials
(see, e.g., B. Lichtin's series of papers such as \cite{Lich88} \cite{Lich93}, 
M. Peter \cite{Pete98}, Essouabri \cite{Esso12}, and so on). 
The contents of those researches are also not treated in the present article.

The series introduced and studied by D. Goldfeld, D. Bump, S. Friedberg, J. Hoffstein et al. 
(see, e.g., \cite{BBF11}) are
also called ``multiple Dirichlet series'', but are different from the series
dsicussed in the present article.

\section{The case of type $*$}

In this section we consider the case of type $*$.    

The simplest situation is when 
$a_j$ ($1\leq j\leq r$) are periodic functions.    
In this case, the sum $\Phi_r^{*}(s_1,\ldots,s_r;a_1,\ldots,a_r)$ can be
easily written as a linear combination of multiple series with numerator $=1$.

The case $a_j=\chi_j$ (Dirichlet characters) was studied by S. Akiyama and
H. Ishikawa \cite{AkIs02}.    They wrote $\Phi_r^{*}(s_1,\ldots,s_r;\chi_1,\ldots,\chi_r)$ 
as a linear combination of multiple zeta-functions of the form
\begin{align}\label{multiple_Hurwitz}
\sum_{m_1=1}^{\infty}\cdots \sum_{m_r=1}^{\infty}
(m_1+\alpha_1)^{-s_1}(m_1+m_2+\alpha_2)^{-s_2}
\cdots (m_1+\cdots+m_r+\alpha_r)^{-s_r}
\end{align}
(where $\alpha_1,\ldots,\alpha_r$ are constants).    The series \eqref{multiple_Hurwitz}
can be treated analogously to the case of \eqref{def_EZ}.
Ishikawa \cite{Ishi01} \cite{Ishi02} studied further, and gave applications to the 
evaluation of certain
multiple character sums.

When the coefficients are not periodic, the study of analytic properties of multiple
series of type $*$ is not easy.
Here we mention a paper of Essouabri, H. Tsumura and the author \cite{EMT11}, in which
the case when the coefficients satisfy a recurrence condition is handled.

Assume the recurrence condition
$$
a_i(m+h)=\sum_{j=0}^{h-1}\lambda_{ji}a_i(m+j)\qquad {\rm for\; all} \quad m\in\mathbb{N}
$$
among coefficients, where $\lambda_{ij}\in\mathbb{C}$.
Denote by $\eta_1,\ldots.\eta_q$ the family of all maps from $\{1,\ldots,r\}$ to
$\{0,\ldots,h-1\}$ (so $q=h^r$), 
and for $\mathbf{m}=(m_1,\ldots,m_r)\in\mathbb{N}^r$, define
$$
A_l(\mathbf{m}):=\prod_{i=1}^r a_i(m_1+\cdots+m_i+\eta_l(i))\qquad (l=1,\ldots,q)
$$
and $A(\mathbf{m}):=(A_1(\mathbf{m}),\ldots,A_q(\mathbf{m})).$
Define the matrix $T_k$ ($1\leq k\leq r$) by
$$
A(\mathbf{m}+h\mathbf{e}_k)=T_k A(\mathbf{m}),
$$
where $\mathbf{e}_k$ is the $k$-th unit vector.

\begin{theorem}
{\rm (\cite{EMT11})}
If $a_i$ is at most of polynomial order, satisfies the above recurrence condition, and
if $1$ is not an eigenvalue of any of the matrices $T_1,\ldots,T_r$, then 
$$
\sum_{m_1=1}^{\infty}\cdots\sum_{m_r=1}^{\infty}
\frac{a_1(m_1)a_2(m_1+m_2)\cdots a_r(m_1+\cdots+m_r)}{m_1^{s_1}(m_1+m_2)^{s_2}
\cdots(m_1+\cdots+m_r)^{s_r}}
$$
can be holomorphically continued to the whole space $\mathbb{C}^r$.
\end{theorem}

As a typical example, a double series with Fibonacci numbers $F_n$ on the numerator
is discussed.
In fact, it was shown that the double series
$$
\phi(s)=\sum_{m=1}^{\infty}\sum_{n=1}^{\infty}\frac{(i\alpha^{-1})^{2m+n}F_{2m+n}
+(i\alpha^{-1})^{m+2n}F_{m+2n}}{m^s(m+n)^s}
$$
($\alpha=(1+\sqrt{5})/2$) can be continued meromorphically to $\mathbb{C}$.
Moreover, 
the evaluation of a special value
$$
\phi(0)=\frac{1}{18}\left(6-\sqrt{5}+(2-3\sqrt{5})i\right)
$$
and a sum formula
$$
\sum_{m=1}^{\infty}\sum_{n=1}^{\infty}\frac{\alpha^{-m}F_m+\alpha^{-n}F_n
-\alpha^{-m-n}F_{m+n}}{m(m+n)^2}
=\sum_{m=1}^{\infty}\frac{\alpha^{-m}F_m}{m^3}
$$
are given.

An important auxiliary tool in \cite{EMT11} is a kind of vectorial zeta-functions.
Special values of those vectorial zeta-functions satisfy ``vectorial sum formulas''.
Such formulas were proved in \cite{EMT11} in the double and triple cases.
The formula in the general case was proposed as a conjecture in \cite{EMT11}, and
proved by S. Yamamoto \cite{Yama15}.

\section{Application of the Mellin-Barnes integral formula}

Hereafter we discuss the case of type $\sh$.
Let $\varphi_k(s)=\sum_{m=1}^{\infty}a_k(m)m^{-s}$ ($1\leq k\leq r$), and we sometimes
write $\Phi_r^{\sh}(s_1,\ldots,s_r;\varphi_1,\ldots,\varphi_r)$ instead of
$\Phi_r^{\sh}(s_1,\ldots,s_r;a_1,\ldots,a_r)$.

The advantage of type $\sh$ is that
this type of series is more suitable to analytic study.   In fact, write
\begin{align*}
&\Phi_r^{\sh}(s_1,\ldots,s_r;\varphi_1,\ldots,\varphi_r)\\
&=\sum_{m_1=1}^{\infty}\cdots\sum_{m_r=1}^{\infty}
\frac{a_1(m_1)a_2(m_2)\cdots a_r(m_r)}{m_1^{s_1}(m_1+m_2)^{s_2}
\cdots(m_1+\cdots+m_{r-1})^{s_{r-1}}}\\
&\times
(m_1+\cdots+m_{r-1})^{-s_r}\left(1+\frac{m_r}{m_1+\cdots+m_{r-1}}\right)^{-s_r},
\end{align*}
and to the last factor
apply the classical Mellin-Barnes integral formula
\begin{align*}
(1+\lambda)^{-s}=\frac{1}{2\pi i}\int_{c-i\infty}^{c+i\infty}
\frac{\Gamma(s+z)\Gamma(-z)}{\Gamma(s)}\lambda^z dz,
\end{align*}
where $s,\lambda\in\mathbb{C}$, $\Re s>0$, $\lambda\neq 0$, $|\arg\lambda|<\pi$, and
$-\Re s<c<0$.
We get
\begin{align*}
&\left(1+\frac{m_r}{m_1+\cdots+m_{r-1}}\right)^{-s_r}\\
&=\frac{1}{2\pi i}\int_{c-i\infty}^{c+i\infty}
\frac{\Gamma(s_r+z)\Gamma(-z)}{\Gamma(s_r)}\left(\frac{m_r}{m_1+\cdots+m_{r-1}}\right)^z dz,
\end{align*}
and hence
\begin{align*}
&\Phi_r^{\sh}(s_1,\ldots,s_r;\varphi_1,\ldots,\varphi_r)\\
&=\frac{1}{2\pi i}\int_{c-i\infty}^{c+i\infty}
\frac{\Gamma(s_r+z)\Gamma(-z)}{\Gamma(s_r)}\\
&\times\sum_{m_1=1}^{\infty}\cdots\sum_{m_{r-1}=1}^{\infty}
\frac{a_1(m_1)a_2(m_2)\cdots a_r(m_r)}{m_1^{s_1}(m_1+m_2)^{s_2}
\cdots(m_1+\cdots+m_{r-1})^{s_{r-1}+s_r+z}}\\
&\qquad\times
\sum_{m_r=1}^{\infty} a_r(m_r)m_r^z dz\\
&=\frac{1}{2\pi i}\int_{c-i\infty}^{c+i\infty}
\frac{\Gamma(s_r+z)\Gamma(-z)}{\Gamma(s_r)}\\
&\qquad\times\Phi_{r-1}^{\sh}(s_1,\ldots,s_{r-2},s_{r-1}+s_r+z;\varphi_1,\ldots,\varphi_{r-1})
\varphi_r(-z)dz.
\end{align*}
(Here, the important point is that the sum with respect to $m_r$ can be separated.)

Using this expression, we can reduce the study of $\Phi_r$ to that of $\Phi_{r-1}$
and $\varphi_r$,
and to that of $\Phi_{r-2}$ and $\varphi_{r-1}$,..... and finally we obtain

\begin{theorem}\label{th_conti}
{\rm (Matsumoto and Tanigawa \cite{MaTa03})}
Assume that $\varphi_k(s)$ ($1\leq k\leq r$) are convergent absolutely for $s$ with sufficiently 
large real part, continued meromorphically to the whole complex plane, with 
finitely many poles, and are of polynomial order.
Then $\Phi_r^{\sh}(s_1,\ldots,s_r;\varphi_1,\ldots,\varphi_r)$ can be continued meromorphically to 
the whole $\mathbb{C}^r$, and the location of possible singularities can be given
explicitly.    In particular, if all $\varphi_k(s)$ ($1\leq k\leq r$) are entire, then
$\Phi_r^{\sh}(s_1,\ldots,s_r;\varphi_1,\ldots,\varphi_r)$ is also entire.
\end{theorem}

In particular,
$$
\Phi_2^{\sh}(s_1,s_2;1,\varphi_2)=\sum_{m_1=1}^{\infty}\sum_{m_2=1}^{\infty}\frac{a_2(m_2)}
{m_1^{s_1}(m_1+m_2)^{s_2}}
$$
can be continued.     
This double series satisfies a certain 
``functional equation'' which is written in terms of confluent
hypergeometric functions (see \cite{ChMa16} \cite{ChMa17}):

\begin{theorem}\label{th_fe}{\rm (Choie and Matsumoto \cite{ChMa16})}
We have
\begin{align*}
\Phi_2^{\sh}(s_1,s_2;1,\varphi_2)=\frac{\Gamma(1-s_1)\Gamma(s_1+s_2-1)}{\Gamma(s_2)}
\varphi_2(s_1+s_2-1)\\
+\Gamma(1-s_1)\left\{F_{2,+}(1-s_2,1-s_1)+F_{2,-}(1-s_2,1-s_1)\right\},
\end{align*}
where
$$
F_{2,\pm}(s_1,s_2)=\sum_{l\geq 1}A_{s_1+s_2-1}(l)\Psi(s_2,s_1+s_2;\pm 2\pi il),
$$
$A_c(l)=\sum_{n|l}n^c a_2(n)$ and $\Psi$ is the confluent hypergeometric function
defined by
$$\psi(a,b;x)=\frac{1}{\Gamma(a)}\int_0^{e^{i\phi}\infty}e^{-xy}y^{a-1}(y+1)^{b-a-1}dy.
$$
\end{theorem}

Furthermore, when $\varphi_2$ is an automorphic $L$-function, then by using the modular
relation, a different type of functional equation can also be proved.

Why we may call Theorem \ref{th_fe} a "functional equation"?
There are mainly two reasons.    First, it can be compared with
the functional equation of the Hurwitz zeta-function 
$\zeta(s,\alpha)=\sum_{n=0}^{\infty}(n+\alpha)^{-s}$ ($0<\alpha\leq 1$):
\begin{align*}
\zeta(1-s,\alpha)=\frac{\Gamma(s)}{(2\pi)^s}\{ e^{\pi is/2}\phi(s,-\alpha)
+e^{-\pi is/2}\phi(s,\alpha)\}
\end{align*}
with $\phi(s,\alpha)=\sum_{n=1}^{\infty}e^{2\pi in\alpha}n^{-s}$.

Secondly, when $a_2(n)\equiv 1$, a symmetric form of the functional equation
(that is, a relation which connects $\Phi_2^{\sh}(s_1,s_2)$ with $\Phi_2^{\sh}(1-s_2,1-s_1)$) 
can be deduced from
Theorem \ref{th_fe} on hyperplanes $s_1+s_2=2k+1$ ($k\in \mathbb{Z}\setminus\{0\}$).

This second point was already observed in a paper of Y. Komori, Tsumura and the author 
\cite{KMT10}.
We note that a generalization of the results in \cite{KMT10} to the case of double 
$L$-functions twisted by Dirichlet characters was discussed in \cite{KMT11}.

In Theorem \ref{th_conti}, the assumption that each $\varphi_k$ has only finitely
many poles is important.    The analytic behavior of $\Phi_r^{\sh}$ may become different
when some of $\varphi_k$ has infinitely many poles.   In the next section we
will discuss such cases.

\section{The case when some $\varphi_k$ has infinitely many poles}

Typical examples of Dirichlet series with arithmetical coefficients which have
infinitely many poles are
$$
\sum_{m=1}^{\infty}\frac{\Lambda(m)}{m^s}=-\frac{\zeta'(s)}{\zeta(s)},
\quad
\sum_{m=1}^{\infty}\frac{\mu(m)}{m^s}=\frac{1}{\zeta(s)},
$$
where
$\Lambda(m)$ is the von Mangoldt function and $\mu(m)$ is the M{\"o}bius function.
Both of the above series have infinitely many poles at the zeros of $\zeta(s)$.

We now consider the behavior of $\Phi_r^{\sh}$, where some of associated $\varphi_k$
has infinitely many poles (as in the above examples).
Let 
$$
N(\Phi_r^{\sh})=\#\{k\;|\; 1\leq k\leq r, \;\varphi_k {\rm \;has \;infinitely \;
many \;poles}\}.
$$.
A recent article of A. Nawashiro, Tsumura and the author \cite{MNT} studied several examples which satisfy 
$N(\Phi_2^{\sh})=1$.    For example,

\begin{theorem}
{\rm (\cite{MNT})}
The double series
$$
\sum_{m_1=1}^{\infty}\sum_{m_2=1}^{\infty}\frac{\Lambda(m_2)}{m_1^{s_1}(m_1+m_2)^{s_2}}
$$
can be continued meromorphiocallty to the whole $\mathbb{C}^2$, and the location of
possible singularities can be described explicitly.
\end{theorem}

Next, take $\varphi_1(s)=\sum_{m=1}^{\infty}a_1(m)m^{-s}$ which has only finitely many
ploes, and consider the convolution of $a_1$ and $\mu$:
$$
\widetilde{a}(m)=\sum_{d|m}a_1\left(\frac{m}{d}\right)\mu(d).
$$
In other words, 
$$ \frac{\varphi_1(s)}{\zeta(s)}=\sum_{m=1}^{\infty}\widetilde{a}(m)m^{-s}.$$
\begin{assumption}\label{Ass}
All non-trivial zeros of $\zeta(s)$ are simple and 
$$
\qquad\frac{1}{\zeta'(\rho_n)}=O(|\rho_n|^B)
\qquad (\rho_n: n{\mbox{-th zero}}, \; B>0)
$$
holds.
(This is a quantitative version of the well-known "simplicity conjecture" for the zeros
of $\zeta(s)$, and
consistent with the Gonek-Hejhal conjecture \cite{Gone89} \cite{Hejh89}.)
\end{assumption}

\begin{theorem}
{\rm (\cite{MNT})}
Under Assumption \ref{Ass}, the double series
$$
\sum_{m_1=1}^{\infty}\sum_{m_2=1}^{\infty}\frac{\widetilde{a}(m_2)}{m_1^{s_1}(m_1+m_2)^{s_2}}
$$
can be continued meromorphiocallty to the whole $\mathbb{C}^2$, and the location of
possible singularities can be described explicitly.
\end{theorem}

The paper \cite{MNT} only considers the double zeta case, but a generalization of
\cite{MNT} to the general multiple case was treated by Rei Kawashima \cite{Kawa19}.
She also discussed the case when $\Lambda$ is replaced by the Liouville function $\lambda$.

\bigskip

The above two theorems show that when $N(\Phi_2^{\sh})=1$, the analytic behavior of
$\Phi_2^{\sh}$ is not so different from the case when $N(\Phi_2^{\sh})=0$.
However, if $N(\Phi_2^{\sh})=2$,  or more generally if $N(\Phi_r^{\sh})\geq 2$, 
the analytic behavior of
the multiple series is totally different.

Consider
$$
F_2(s)=\sum_{k=1}^{\infty}\sum_{l=1}^{\infty}\frac{\Lambda(k)\Lambda(l)}{(k+l)^s},
$$
which is equal to $\Phi_2(0,s;-\zeta'/\zeta,-\zeta'/\zeta)$, so
both $\varphi_1$ and $\varphi_2$ have infinitely many poles.    
We can rewrite
$$
F_2(s)=\sum_{m=1}^{\infty}\frac{G_2(m)}{m^s},
\qquad G_2(m)=\sum_{k+l=m}\Lambda(k)\Lambda(l).
$$
It is to be noted that this $G_2(m)$ is the counting function of the classical
Goldbach problem:
$$
G_2(m)=\sum_{r_1\geq 1, r_2\geq 1 \atop{p_1^{r_1}+p_2^{r_2}=m}}\log p_1 \log p_2
=\sum_{p_1+p_2=m}\log p_1 \log p_2+({\mbox{error}}),
$$
where $p_1,p_2$ denote primes.
Using the Mellin-Barnes formula, we have
$$
F_2(s)=\frac{1}{2\pi i}\int_{c-i\infty}^{c+i\infty}\frac{\Gamma(s+z)\Gamma(-z)}
{\Gamma(s)}\frac{\zeta'}{\zeta}(s+z)\frac{\zeta'}{\zeta}(-z)dz.
$$
Shifting the path of integration suitably, we find that there are poles of $F_2(s)$ at:
$$
s=2, \quad s=\rho+1,\quad  s=\rho+\rho'
$$
(where $\rho$, $\rho'$ denotes the non-trivial zeros of $\zeta(s)$).

Now assume the RH (Riemann Hypothesis).    Then $\Re(\rho+\rho')=1$, and we
can show that $\rho+\rho'$ are dense on the line $\Re s=1$.
Therefore $\Re s=1$ seems a kind of barrier if we want to continue $F_2(s)$
meromorphically.

It is believed that $\gamma=\Im\rho>0$ are linearly independent over $\mathbb{Q}$
(the linear independence conjecture, LIC).

\begin{theorem} If RH and LIC are true, then $\Re s=1$ is the natural boundary of
$F_2(s)$.
\end{theorem}

This was first proved by S. Egami and the author \cite{EgMa07} under the RH and 
a stronger quantitative version of LIC, and then in the above form by G. Bhowmik and J.-C. Schlage-Puchta \cite{BhSc11}.
Moreover, in \cite{BhSc11} it is also shown that
$$
F_r(s)=\sum_{k_1=1}^{\infty}\cdots\sum_{k_r=1}^{\infty}\frac{\Lambda(k_1)\cdots
\Lambda(k_r)}{(k_1+\cdots+k_r)^s}\qquad (r\geq 3)
$$
has the natural boundary $\Re s=r-1$ if RH is true and $\Re s=1$ is the natural boundary of
$F_2(s)$.

A generalization to the case with congruence conditions has been also studied.    Let
$$
G_2(m;q,a,b):=\sum_{\substack{k+l=m \\k\equiv a,l\equiv b\; ({\rm mod}q)}}\Lambda(k)\Lambda(l),
$$
where $a,b,q$ are positive integers with $(ab,q)=1$, 
and define
the associated Dirichlet series by
$$
F_2(s;q,a,b):=\sum_{m=1}^{\infty}\frac{G_2(m;q,a,b)}{m^s}.
$$
The behavior of $F_2(s;q,a,b)$ was first treated by F. R{\"u}ppel \cite{Rupp12}, and then 
further studied by Y. Suzuki \cite{Suzu17}.

Let
$$
S(x;q,a,b):=\sum_{m\leq x}G_2(m;q,a,b).
$$
We can reduce the study of $S(x;q,a,b)$, via its associated Dirichlet series
$F_2(s;q,a,b)$, to that of the behavior (especially the distribution of zeros) of
Dirichlet $L$-functions.
Thereby we can establish the connection between Goldbach generating functions
and the zeros of $L$-functions (especially the generalized Riemann hypothesis, GRH).
The existence of this connection was
first suggested by A. Granville \cite{Gran07} in the Riemann zeta case,
and then fully developed in Bhowmik et al.
\cite{BhHa20} \cite{BHMS19} \cite{BhRu18}.
We conclude this article with the statement of some theorems proved in those papers.

Let $B_{\chi}=\sup\{\Re\rho_{\chi}\}$ and
$B_q=\sup\{B_{\chi}\;|\;\chi \;({\mbox{mod}} \;q)\}$.

\begin{theorem}
{\rm (\cite{BHMS19})}
For any $\delta>0$, We have
$$
S(x;q,a,b)=\frac{x^2}{2\varphi(q)^2}+O(x^{1+B_q^*}),
$$
where $B_q^*=\min\{B_q,1-\eta\}$ with
$$
\eta=\frac{c(\delta)}{\max\{q^{\delta},(\log x)^{2/3}(\log\log x)^{1/3}\}}
$$
{\rm (}where $c(\delta)$ is a small positive constant{\rm )}.
\end{theorem}

\begin{remark}
In particular, if we assume GRH, then $B_q^*=B_q=1/2$ and
$$
S(x;q,a,b)=\frac{x^2}{2\varphi(q)^2}+O(x^{3/2}).
$$
\end{remark}

Now recall a well-known conjecture (DZC): Any two distinct Dirichlet 
$L$-functions (mod $q$) do not have a common 
non-trivial zero (except for a possible multiple zero at $s=1/2$).

\begin{theorem}
{\rm ( \cite{BHMS19}, \cite{BhRu18})}
Assume that the DZC is true, and $\chi(a)+\chi(b)\neq 0$ for all $\chi$ {\rm(}mod $q${\rm)}. 
If the asymptotic formula
$$
S(x;q,a,b)=\frac{x^2}{2\varphi(q)^2}+O(x^{1+d+\varepsilon}) \quad (1/2\leq d<1)
$$
holds for any $\varepsilon>0$, then either $B_q\leq d$ or $B_q=1$.   Moreover
we can remove the possibility of $B_q=1$ when $a=b$.
\end{theorem}

\begin{remark}
In particular, if $a=b$ and the above formula holds with
the error $O(x^{3/2+\varepsilon})$ (that is, $d=1/2$), then $B_q=1/2$, that is, 
the GRH holds.
\end{remark}


\end{document}